\documentclass[journal]{IEEEtran}
\usepackage{graphicx}
\ifCLASSINFOpdf
\else
\fi
\usepackage{wrapfig}

\usepackage{amsthm, amssymb, amsmath, latexsym, amsfonts, mathcomp}
\theoremstyle{definition}
\newtheorem{remark}{Remark}

\usepackage{fixltx2e}
\begin{document}
%
\title{Polar Wavelets in Space}
%
%
%

\author{Christian Lessig
\thanks{Christian Lessig is with the Institute of Simulation and Graphics at Otto-von-Guericke-Universit{\"a}t Magdeburg, Germany.}
}

%
%

\markboth{IEEE Signal Processing Letters}%
{Lessig: Polar Wavelets in Space}
%



\maketitle

\begin{abstract}
Recent work introduced a unified framework for steerable and directional wavelets in two and three dimensions that ensures many desirable properties, such as a multi-scale structure, fast transforms, and a flexible angular localization.
We show that, for an appropriate choice for the radial window function, these wavelets also have closed form expressions for, among other things, the spatial representation,  the filter taps for the fast transform, and the frame representation of the Laplace operator.
The numerical practicality and benefits of our work are demonstrated using signal estimation from non-uniform, point-wise samples, as required for example in ray tracing, and for reconstructing a signal over a lower-dimensional sub-manifold, with applications for instance in medical imaging.
\end{abstract}


%
\IEEEpeerreviewmaketitle

\section{Introduction}
\label{sec:introduction}

Steerable wavelets have a long history, going back to work by Freeman and Adelson~\cite{Freeman1991} and Perona~\cite{Perona1991}.
Recently, a systematic framework for their construction in two and more dimensions has been introduced by Unser and co-workers~\cite{Unser2011,Chenouard2012,Unser2013,Ward2014}, subsuming also later directional wavelets such as curvelets~\cite{Candes2005a,Candes2005b} and related constructions~\cite{Do2005a,Labate2005}.
The work by Unser et al. builds on isotropic multi-resolution analyses~\cite{Portilla2000,Papadakis2003,Azencott2009,Romero2009a} that are defined through a radial window $\hat{h}(\vert \xi \vert)$ (see inset figure).
\setlength{\columnsep}{4pt}
\begin{wrapfigure}[10]{r}[0pt]{0.5\columnwidth}
  \centering
  \vspace{-11pt}
  \includegraphics[width=0.5\columnwidth]{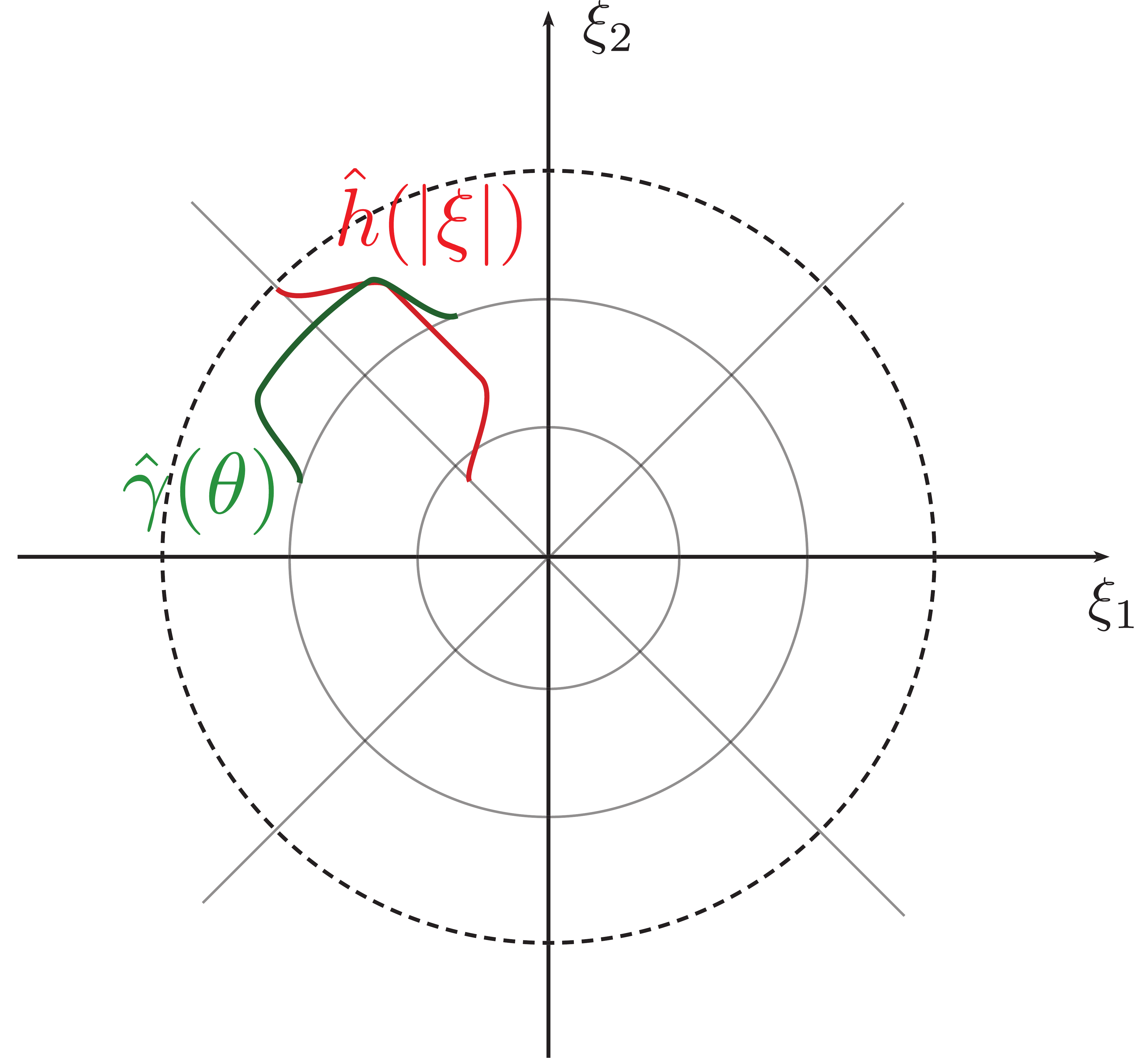}
  \vspace{-10pt}
\end{wrapfigure}
Steerability and directionality are introduced through an angular localization window $\hat{\gamma}( \theta_{\xi})$.
The wavelets are then defined as $\hat{\psi}(\xi) = \hat{\gamma}( \theta_{\xi} ) \, \hat{h}(\vert \xi \vert)$ and since their construction relies on polar (and spherical) coordinates, we will refer to them as polar wavelets, or short polarlets.

Under suitable admissibility conditions for the window functions (cf. Sec.~\ref{sec:construction}) polarlets generate a Parseval tight frame, have a multi-resolution structure, and fast transforms.
The present letter extends the useful properties of polar wavelets even further.
We show that, with the Portilla-Simoncelli radial window~\cite{Portilla2000}, the spatial representation of the wavelets and the filter taps for the fast transforms have closed form expressions.
Furthermore, also calculations like the Galerkin projection of the Laplace operator, which is required, for example, for the solution of Poisson problems, can be performed analytically.
Although the closed form expressions are complicated, compared to the numerical quadrature which would be required otherwise, they enable theoretical insights and more easily provide the accuracy that is required.
In the three dimensional setting, we also provide more concrete and practical admissibility conditions for the angular localization window and a proof that, in our opinion, is simpler than existing ones in the literature~\cite{Chenouard2012,Ward2014}.

To demonstrate the relevance of our work we provide two proof-of-concept applications that are facilitated by our results: the estimation of an image signal from pointwise samples, as used in ray tracing-based image generation, and the reconstruction of a signal on a sub-manifold, which is of importance for example in medical imaging.

\begin{figure}[t]
  \includegraphics[trim={0, 0, 0, 0}, clip, width=0.48\columnwidth]{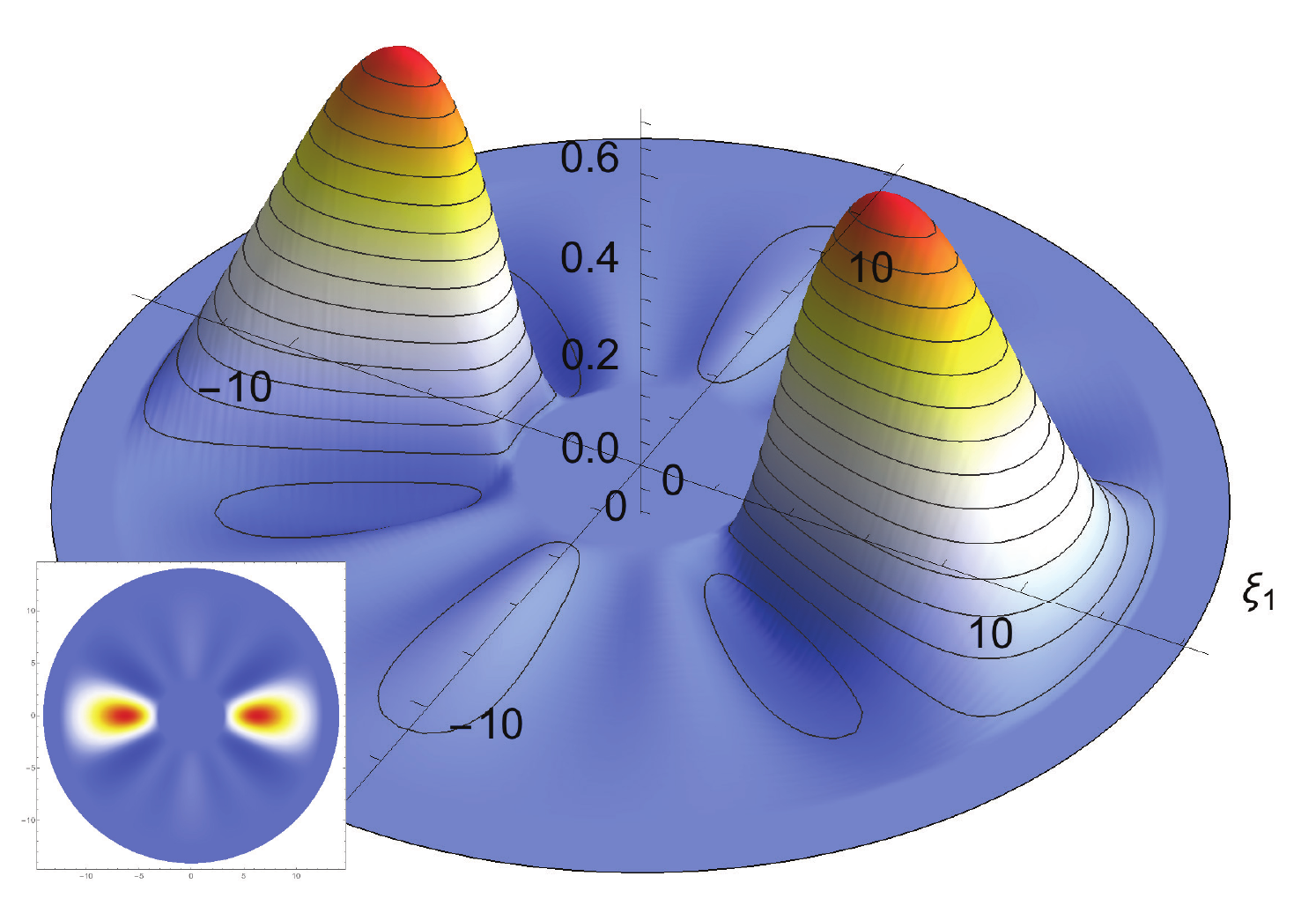}
  \includegraphics[trim={0, 0, 0, 0}, clip, width=0.48\columnwidth]{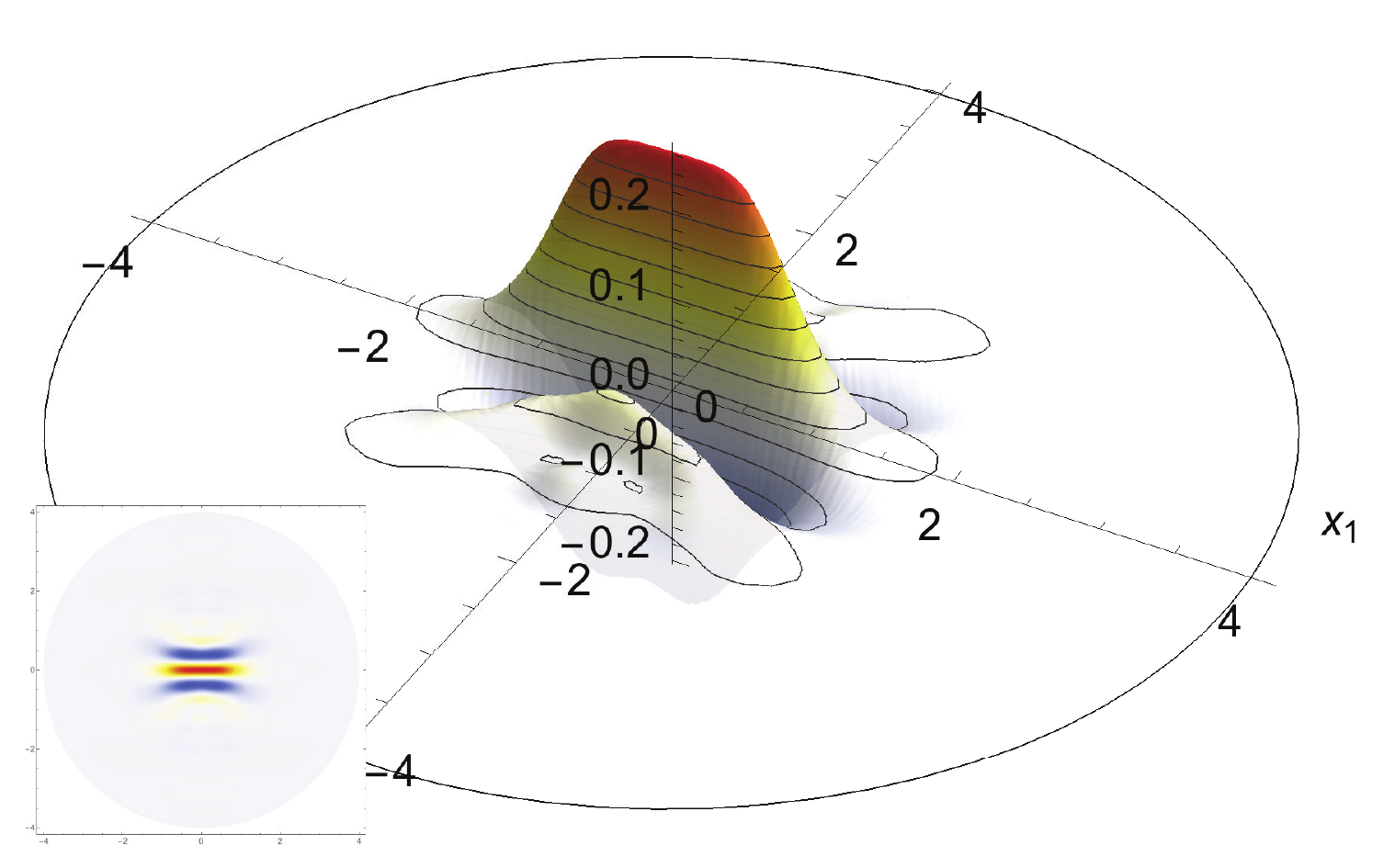}
  \caption{Directional wavelet in two dimensions in the frequency (left) and spatial (right) domains. For the radial window $\hat{h}(\vert \xi \vert)$ we use those proposed by Portilla and Simoncelli~\cite{Portilla2000} and the angular one is a modification of the orthonormal wavelets for $S^1$ by Walter and Cai~\cite{Walter1999}. The wavelets have closed form expressions in both the space and frequency domains.}
  \label{fig:psi:2d:j1}
\end{figure}

%
%

\section{Polar Wavelets}
\label{sec:construction}

\subsection{Construction in Two Dimensions}
\label{sec:construction:2d}

Polar wavlets in two dimensions are defined by~\cite{Unser2013}
\begin{align}
  \label{eq:2D:psi:hat}
  \hat{\psi}(\xi)
  = \hat{\gamma}(\theta_{\xi}) \, \hat{h}(\vert \xi \vert)
  = \left( \sum_{n} \beta_{j,n} \, e^{i n \theta_{\xi}} \right) \, \hat{h}(\vert \xi \vert)
\end{align}
and they generate the spatial wavelet functions
\begin{align}
  \label{eq:2D:psi}
  \psi_{jkt}(x) \equiv \frac{2^j}{2\pi} \psi \big( R_{jt} \, 2^{j} x - k \big)
\end{align}
where $\psi(x)$ is the inverse Fourier transform of Eq.~\ref{eq:2D:psi:hat}, $k \in \mathbb{Z}^2$ is a translation, and $R_{jt}$ an equal angle rotation matrix by $2\pi / M_j$, with $M_j$ being the number of different orientations at each $k \in \mathbb{Z}^2$.
Unser and Chenouard~\cite{Unser2013} showed that the functions in Eq.~\ref{eq:2D:psi} form a Parseval tight frame when the radial window $\hat{h}(\vert \xi \vert)$ satisfies the Cald{\`e}ron admissibility condition
\begin{align}
  \label{eq:admissibility:calderon:radial}
  \sum_{j \in \mathbb{Z}} \big\vert \hat{h}\big( 2^{-j} \vert \xi \vert \big) \big\vert^2 = 1 \quad , \quad \forall \xi \in \mathbb{R}_{\xi}^2
\end{align}
and the matrix $U$, formed by all Fourier series coefficients $\beta_{j,n}^t = R_{jt}^T \, \beta_{j,n} = e^{i n t (2\pi / M_j)} \, \beta_{j,n}$, fulfills $U^H U = D$ with the diagonal matrix $D$ having $\mathrm{tr}( D ) = 1$ (as usual, in practice one uses also scaling functions; the following discussion applies with obvious modifications also to these).
Unser and Chenouard~\cite{Unser2013} establish the result using the higher order Riesz transform.
With a slight generalization of the admissibility condition in Eq.~\ref{eq:admissibility:calderon:radial}, the result can also be obtained using a direct computation, see the supplementary material.

\begin{remark}
Intuition for the condition $U^H U = D$ can be obtained through the special case $U^H U = \mathrm{Id} / (2N+1)$ which, up to a rescaling, is satisfied when the $\beta_{j,n}^t$ generate a tight frame for the space generated by the Fourier series up to band $N$.
In our examples we use as angular localization window hence a variation of the orthonormal wavelets for $S^1$ by Walter and Cai~\cite{Walter1999}, which are simple yet general.
Our modification symmetrizes the scaling functions by a rotation by $\pi$, so that one obtains real-valued $\psi_{jkt}(x)$, see Fig.~\ref{fig:psi:2d:j1} for an example.
\end{remark}

\paragraph*{Spatial representation}
The compute the inverse Fourier transform of Eq.~\ref{eq:2D:psi:hat} it is convenient to work with the Jacobi-Anger formula,
\begin{align}
  \label{eq:jacobi_anger}
  e^{i \langle x , \xi \rangle} = \sum_{m \in \mathbb{Z}} i^m \, e^{i m (\theta_{x} - \theta_{\xi})} J_m(\vert \xi \vert \, \vert x \vert ) ,
\end{align}
which expresses the complex exponential in polar coordinates.
With it, and by evaluating the integral in the inverse Fourier transform in polar coordinates, one obtains
\begin{subequations}
  \label{eq:2D:psi:hat:inverse}
\begin{align}
  \psi(x)
  &= \frac{1}{2 \pi} \int_{\mathbb{R}_{\xi}^2} \hat{\psi}(\xi) \, e^{i \langle \xi , x \rangle} \, d\xi
  \\[3pt]
  &= \frac{1}{2 \pi} \sum_{n} \sum_{m \in \mathbb{Z}} i^m \beta_{j,n} \, e^{i n \theta_x} \int_{S_{\theta_{\xi}}^1} \, e^{i n \theta_{\xi}} \, e^{-i m \theta_{\xi}}  d\theta_{\xi}
  \\
  & \quad \quad \quad \quad \quad \quad \quad \times \int_{\mathbb{R}_{\vert \xi \vert}^+} \hat{h}(\vert \xi \vert) \, J_m(\vert \xi \vert \, \vert x \vert ) \, d\vert \xi \vert
  \nonumber
  \\[3pt]
  \label{eq:2D:psi:hat:inverse:3}
  &= \sum_{n} i^n \beta_{j,n} \, e^{i n \theta_x}
  \underbrace{\int_{\mathbb{R}_{\vert \xi \vert}^+} \hat{h}(\vert \xi \vert) \, J_n(\vert \xi \vert \, \vert x \vert ) \, d\vert \xi \vert}_{h_n(\vert x \vert)} .
\end{align}
\end{subequations}
The angular localization window is hence invariant under the inverse Fourier transform, up to the factor of $i^n = e^{i n \pi / 2}$ which implements a rotation by $\pi / 2$, see again Fig.~\ref{fig:psi:2d:j1}.

%

\begin{figure}
  \includegraphics[width=\columnwidth]{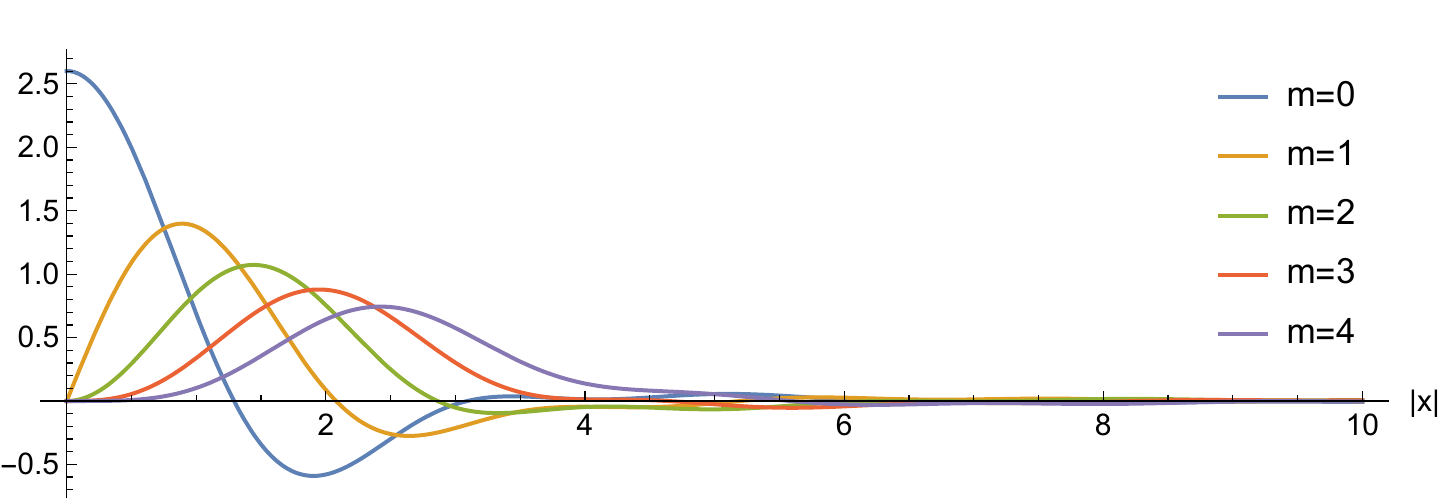}
  \caption{Radial window functions $h_m(\vert x \vert)$ in the spatial domain whose closed form expression is given in Eq.~\ref{eq:h:explicit}.}
  \label{fig:polarlets:hk}
\end{figure}

Eq.~\ref{eq:2D:psi:hat:inverse:3} appeared before in~\cite{Unser2013}.
However, to obtain an explicit expression for the inverse Fourier transform we also have to evaluate the Hankel transform that is given by the radial integral in Eq.~\ref{eq:2D:psi:hat:inverse:3}.
To our knowledge, it has not been observed that for the Portilla-Simoncelli window,
\begin{align}
  \label{eq:def:h}
  \hat{h}(r) =
  \left\{
    \begin{array}{cc}
      \cos{\left( \frac{\pi}{2} \log_2\left(\frac{2 \vert \xi \vert}{\pi} \right)\right)}
      & \frac{\pi}{4} < \vert \xi \vert < \pi
      \\[5pt]
      0 & \textrm{otherwise}
    \end{array}
  \right.
\end{align}
it has a closed form expression,
\begin{subequations}
  \label{eq:h:explicit}
\begin{align}
  h_n(r) &= -i \, 8^{-2 - n} \pi^{2 + n} \, r^n
  \\[4pt]
  &\Bigg(
  \Gamma(c_2^-) \Big( 4^{n+2} \, \tilde{G}_{2,4}^{4,-}(r) + \, \tilde{G}_{2,4}^{64,-}(r) \Big) -
  \nonumber
  \\
  & \ \ \ \Gamma(c_2^+) \Big( 4^{n+2} \, \tilde{G}_{2,4}^{4,+}(r) + \, \tilde{G}_{2,4}^{64,+}(r)
  \Big) \Bigg)
  \nonumber
 \\[3pt]
 & c_b^{\pm} = \frac{1}{2} \left( b + n \pm \frac{i \pi }{\log (4)}\right)
 \\[3pt]
 & \tilde{G}_{a,b}^{c,\pm}(r) = \, _1\tilde{F}_2\left( c_a^{\pm} ; n+1, c_b^{\pm} ; -\frac{1}{c} \pi ^2 r^2\right)
\end{align}
\end{subequations}
where $_1\tilde{F}_2$ is the modified hypergeometric function, see Fig.~\ref{fig:polarlets:hk} for plots.
While, without doubt, quite complicated, Eq.~\ref{eq:h:explicit} can be evaluated to arbitrary precision and it can be analyzed using the extensive results on hypergeometric functions that are available.
Next to the Portilla-Simoncelli window in Eq.~\ref{eq:def:h} we also considered the other radial functions listed in~\cite{Unser2011}.
To the best of our knowledge, however, none of them has a closed form solutions for the Hankel transform.

\paragraph*{Fast transform}
The filter taps for the fast transform are
\begin{subequations}
\begin{align}
  \alpha_{j,k} &= \big\langle \phi_{j,0}(x) , \phi_{j+1,k}(x) \big\rangle
  \\[2pt]
  \beta_{j,k,t} &= \big\langle \psi_{j,0,t}(x) , \phi_{j+1,k}(x) \big\rangle .
\end{align}
\end{subequations}
Using Parseval's theorem and computing the inner product in the Fourier domain, with the translation term expanded using the Jacobi-Anger formula, one obtains for $\beta_{j,k,t}$ that
\begin{subequations}
  \label{eq:fast_transform}
\begin{align}
  \beta_{j,k,t} = \sum_{n} \, i^n \, \beta_{j,n}^t \, e^{i n \theta_k} \, B_n(\vert k \vert) .
\end{align}
with the radial term given by
\begin{align}
  B_{n}(\vert k \vert) = \! \! \int_{\mathbb{R}_{\vert \xi \vert}^+} \! \! \hat{h} \big(\vert 2^{-j} \xi \vert\big) \hat{g} \big(\vert 2^{-j-1} \xi \vert \big) \, J_n\big(\vert \xi \vert \, \vert 2^{-j-1} k \vert \big) \,  \vert \xi \vert \, d\vert \xi \vert
  \nonumber
\end{align}
\end{subequations}
$B_{n}(\vert k \vert)$ can be evaluated in closed form, see the supplementary material for the expression and a detailed derivation as well as those for $\alpha_{j,k}$.

\paragraph*{Other results}
Closed form solutions can also be obtained for other common calculations.
As an example we consider the Galerkin projection of the Laplace operator $\Delta$.
As shown in the supplementary material, for isotropic wavelets it is
\begin{align}
  v_{sr}
  &= \big\langle  \Delta \psi_s \, , \, \psi_r \big\rangle
  \\
  &= -2\pi \int_{\mathbb{R}_{\vert \xi \vert}^+} \hat{h}( 2^{-j_s} \vert \xi \vert) \, \hat{h}( 2^{-j_r} \vert \xi \vert) \, J_0( \vert \xi \vert \, \vert 2^{-j_r} k_r \vert )  \, \vert \xi \vert^3 \, d\vert \xi \vert
  \nonumber
\end{align}
with the radial integral again having a closed form solution.
Analogous expressions can be derived for non-isotropic $\psi_{j,k,t}(x)$.


\subsection{Construction in Three Dimensions}
\label{sec:construction:3d}

{\setlength{\abovecaptionskip}{0pt}
{\setlength{\belowcaptionskip}{0pt}
\begin{figure}
  \includegraphics[trim={10 180 0 85},clip,width=0.48\columnwidth]{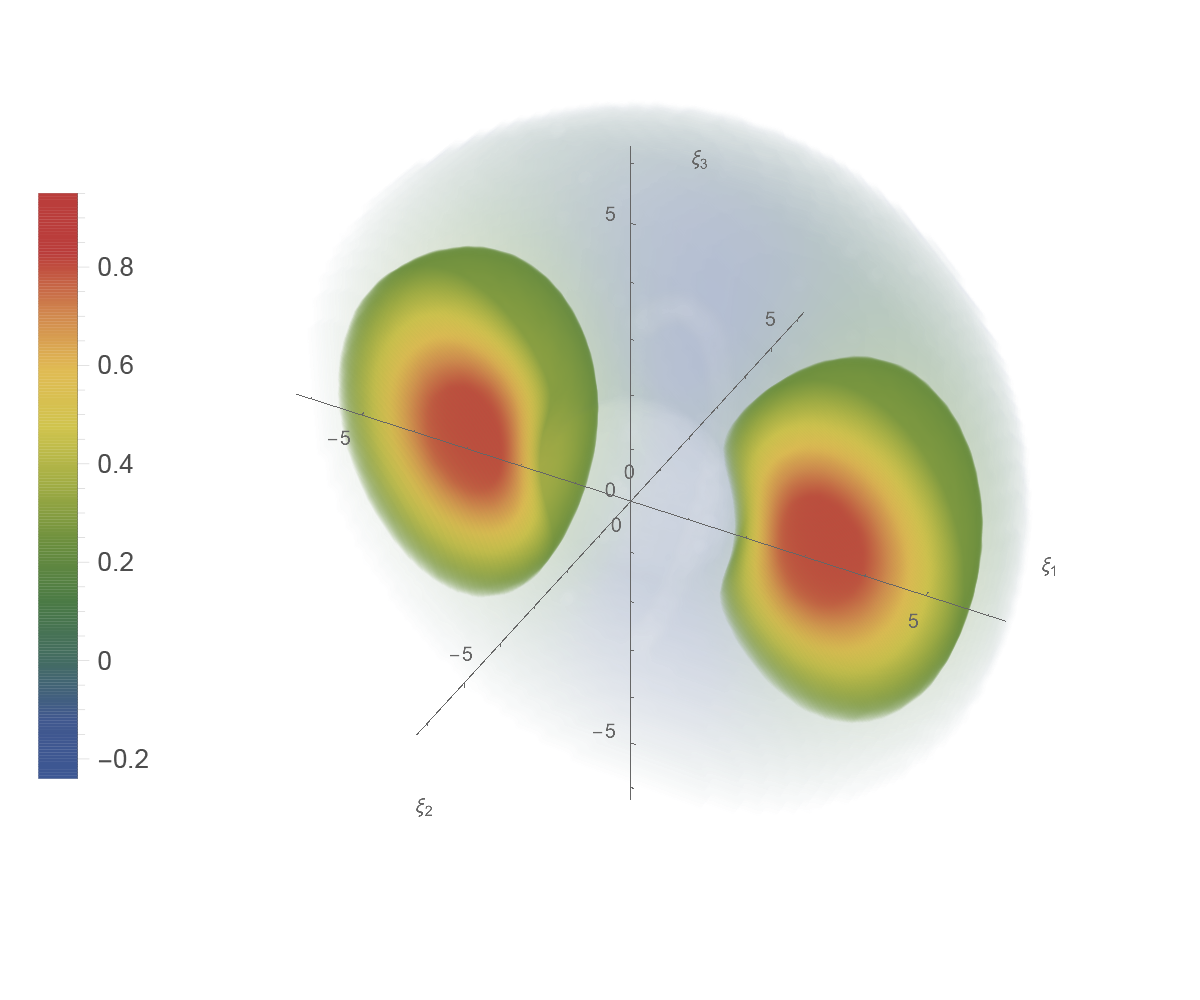}
  \includegraphics[trim={0 180 0 85},clip,width=0.48\columnwidth]{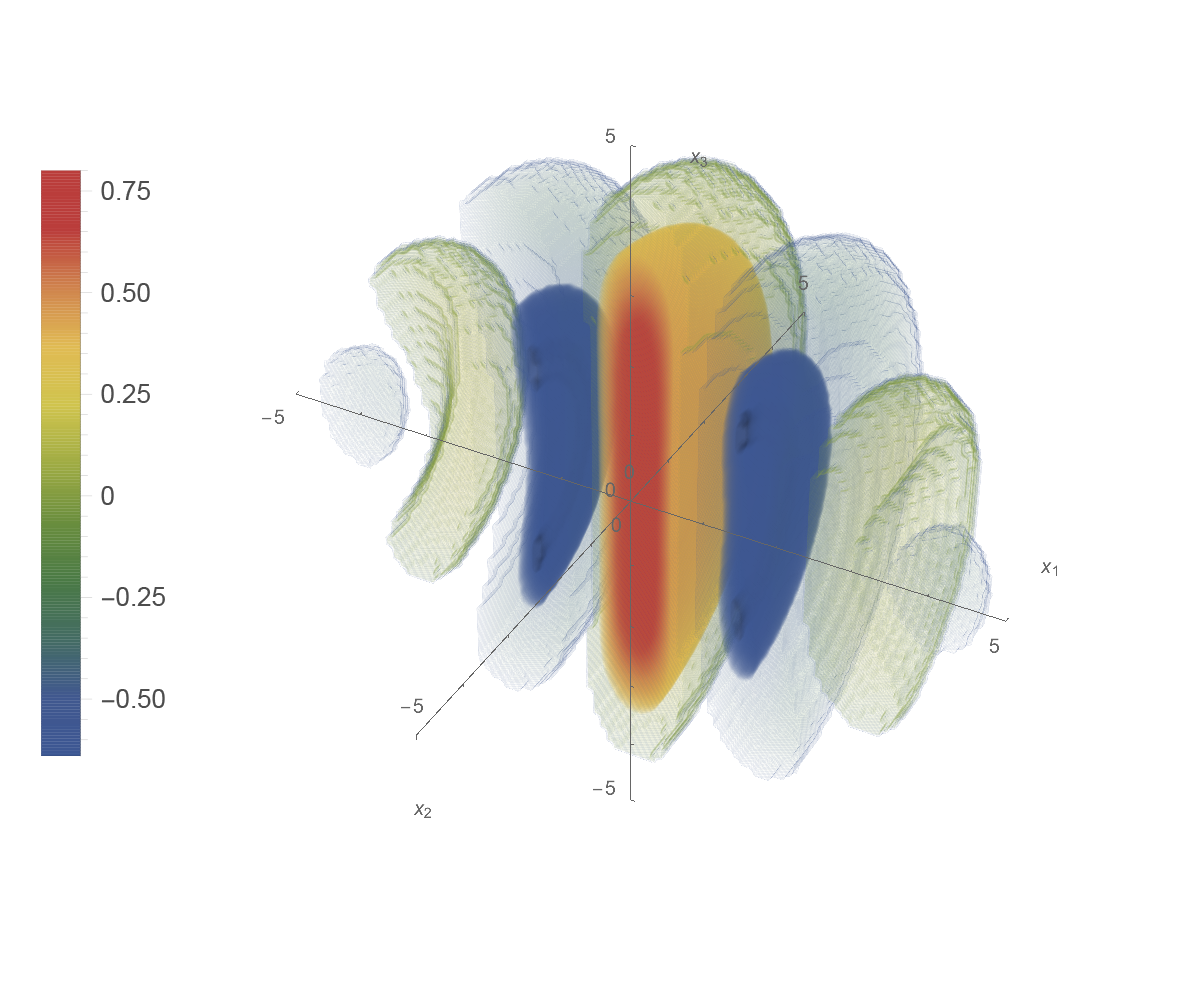}
  \includegraphics[trim={0 0 0 0},clip,width=0.46\columnwidth]{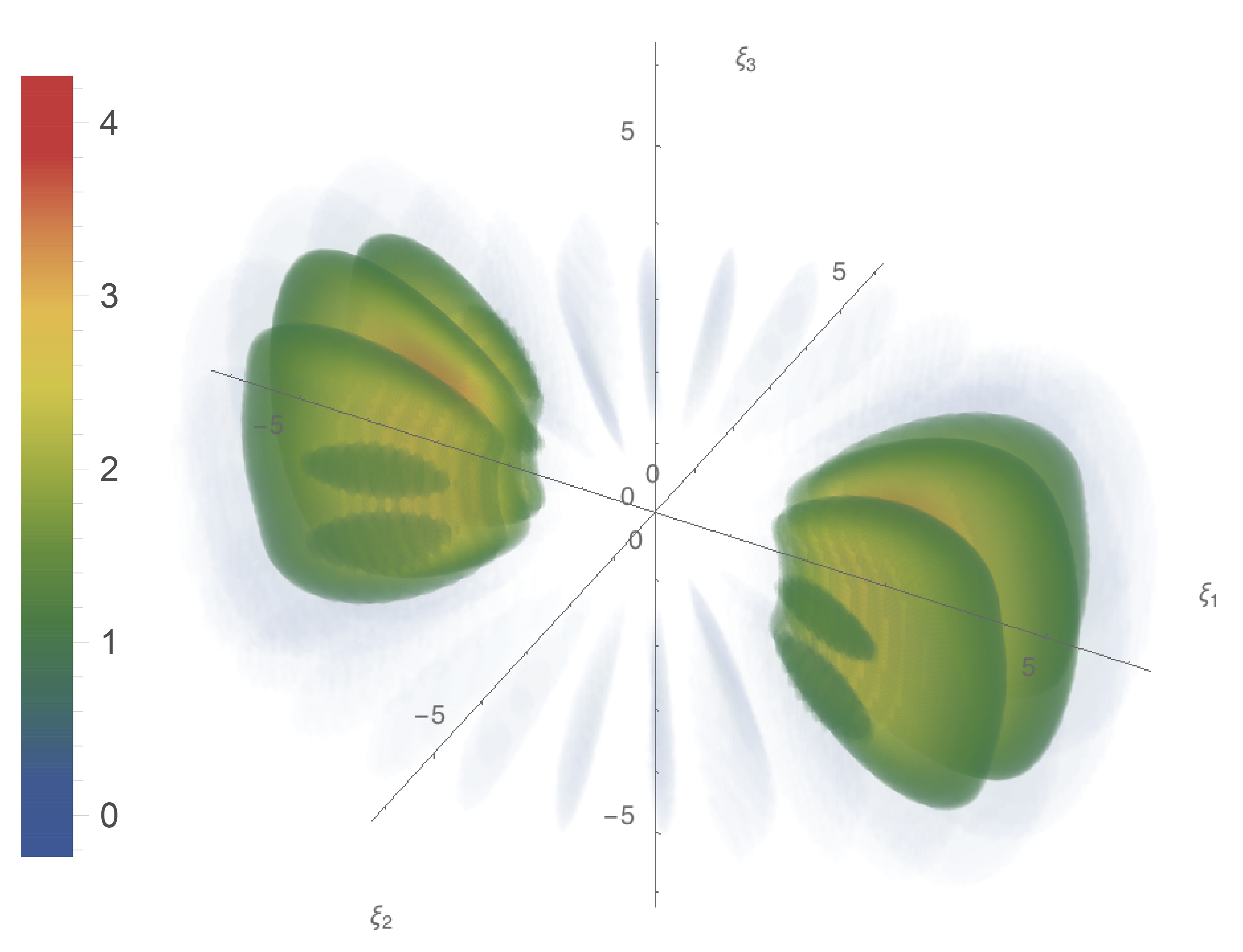}
  \hspace{0.1in}
  \includegraphics[trim={0 0 0 0},clip,width=0.48\columnwidth]{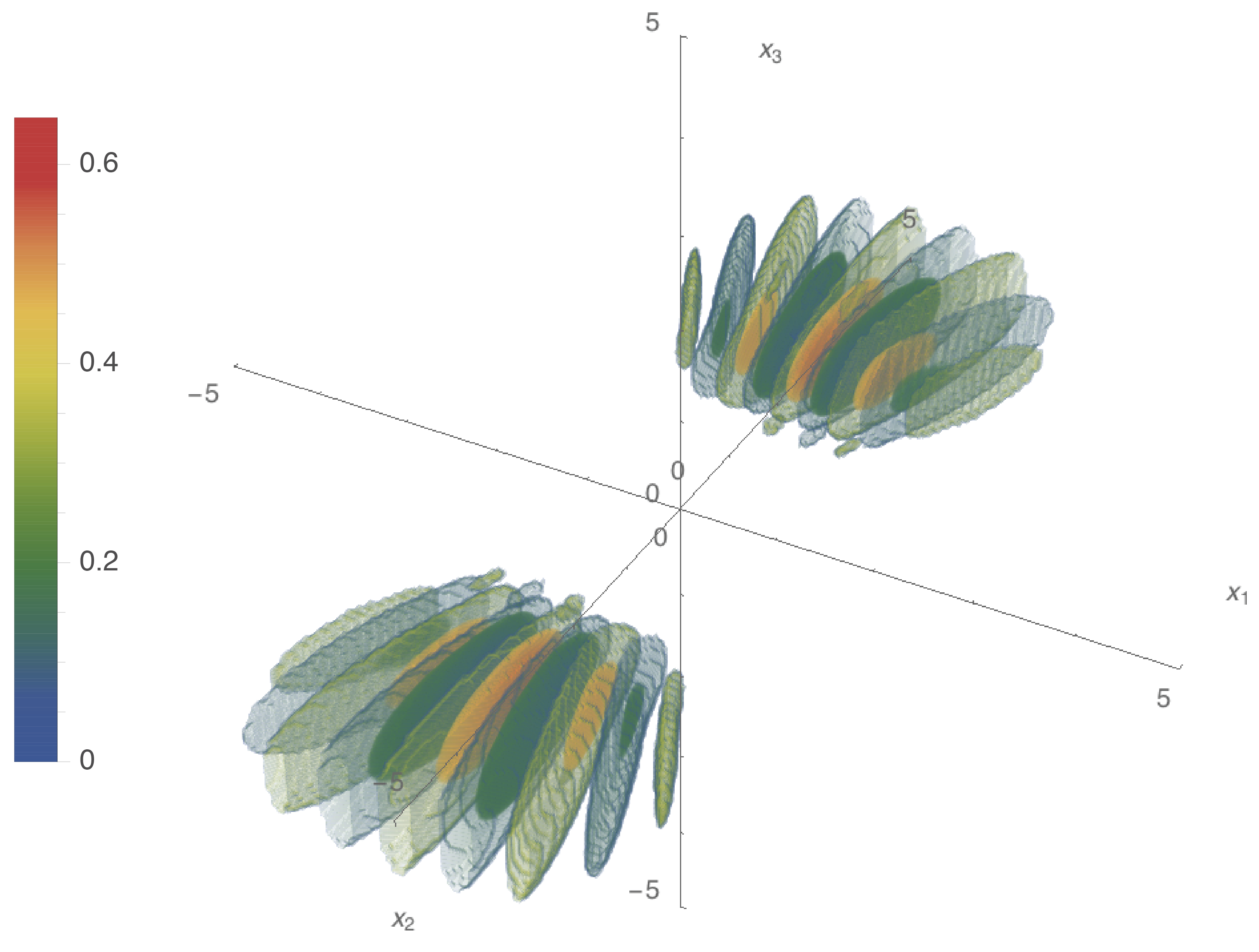}
  \caption{Directional frame functions in three dimensions in the Fourier domain (left) and spatial domain (right) using isotropic (top) and anisotropic (bottom) localization windows. With the isotropic windows from~\cite{Narcowich2006a} and the angular localization windows of Sec.~\ref{sec:construction:2d} they define a Parseval tight frame.}
  \label{fig:psi_3d}
\end{figure}
}

In three dimensions, polar wavelets are given by~\cite{Chenouard2012,Ward2014}
\begin{align}
  \label{eq:psi:hat:3d}
  \hat{\psi}(\xi)
  = \hat{\gamma}(\bar{\xi}) \, \hat{h}(\vert \xi \vert)
  = \left( \sum_{lm} \kappa_{lm}^{tj} \, y_{lm}(\bar{\xi}) \right) \hat{h}( \vert \xi \vert )
\end{align}
where the $y_{lm}(\bar{\xi})$, with $\bar{\xi} = \xi / \vert \xi \vert$, are spherical harmonics~\cite{Freeden1998}.
When the window coefficients $\kappa_{lm}$ only depend on $l$, i.e. $\kappa_{lm} = \bar{\kappa}_l$ as in~\cite{Narcowich2006a,McEwen2016}, then one obtains isotropic wavelets, see Fig.~\ref{fig:psi_3d}, top.
Anisotropic curvelet- and ridgelet-like window functions are obtained using $\kappa_{lm} = \beta_m \, \bar{\kappa}_l$, where the $\beta_m$ are for example again the windows used in the last section.
With such $\kappa_{lm}$ one obtains analogues of the functions discussed by Kutyniok and Petersen~\cite[Sec. 4.3]{Kutyniok2016} that are able to extract curvature information, see Fig.~\ref{fig:psi_3d}, bottom.
Unser, Chenouard and Ward~\cite{Ward2014,Chenouard2012} used again the higher-order Riesz transform to obtain admissibility conditions so that Eq.~\ref{eq:psi:hat:3d} generates a tight frame.
In the supplementary material we provide a direct argument that shows that again Eq.~\ref{eq:admissibility:calderon:radial} has to be satisfied and that the angular windows have to fulfill
\begin{align}
  \label{eq:admissibility:3d:angular}
  \delta_{l,0} \delta_{m,0} = \sum_{t} u_{jt}^T \, G^{lm} \, u_{jt}
\end{align}
where the $u_{jt}$ are the vectors formed by the $\kappa_{lm}^{tj}$.
To our knowledge, Eq.~\ref{eq:admissibility:3d:angular} did not appear before in the literature.

\paragraph*{Spatial representation}
Previous work~\cite[Theorem 2.5 and Theorem 4.2]{Ward2014} established some properties of three-dimensional polarlets in the spatial domain.
However, to our knowledge currently no explicit description of $\psi(x)$ exists.
We show in the following that it can be obtained using a computation similar to Eq.~\ref{eq:2D:psi:hat:inverse}.

In three dimensions, the role of the Jacobi-Anger formula is taken by the Rayleigh formula
\begin{align}
  e^{i \langle \xi , x \rangle} = 4 \pi \sum_{l=0}^{\infty} \sum_{m=-l}^l i^l \, y_{lm}( \bar{x} ) \, y_{lm}^*( \bar{\xi} ) \, j_l( \vert \xi \vert \, \vert x \vert )
\end{align}
where $j_l( \vert \xi \vert \, \vert x \vert )$ is the spherical Bessel function of degree $l$.
By evaluating the integral of the inverse Fourier transform of Eq.~\ref{eq:psi:hat:3d} in spherical coordinates one obtains
\begin{subequations}
  \label{eq:psi:3d:fourier_inverse}
\begin{align}
  \label{eq:psi:3d:fourier_inverse:1}
  \psi(x)
  &= \int_{\mathbb{R}_{\xi}^3} \hat{\psi}(\xi) \, e^{i \langle \xi , x \rangle} \, d\xi
  \\
  &= \sum_{l',m'} i^{l'}  y_{l'm'}( \bar{x} ) \, \sum_{lm} \kappa_{lm}
  \int_{S_{\bar{\xi}}^2} y_{lm}(\bar{\xi}) \, y_{l'm'}^*( \bar{\xi} ) \, d\bar{\xi}
  \nonumber
  \\
  \label{eq:psi:3d:fourier_inverse:2}
  & \quad \quad \quad \quad \quad \times \int_{\mathbb{R}_{\vert \xi \vert}^+} \hat{h}(\vert \xi \vert) \, j_{l'}( \vert \xi \vert \, \vert x \vert ) \vert \xi \vert^2 d\vert \xi \vert .
\end{align}
With the orthonormality of the spherical harmonics, we have
\begin{align}
  \label{eq:psi:3d:fourier_inverse:3}
  \psi(x) = \sum_{l,m} i^l \kappa_{lm} \, y_{lm}( \bar{x} ) \underbrace{\int_{\mathbb{R}_{\vert \xi \vert}^+} \hat{h}(\vert \xi \vert) \, j_{l}( \vert \xi \vert \, \vert x \vert ) \vert \xi \vert^2 d\vert \xi \vert}_{h_l(\vert x \vert)}
\end{align}
\end{subequations}
Analogous to the two-dimensional setting the angular window is invariant under the Fourier transform up to a factor of $i^l$.
It again implements a rotation by $\pi / 2$, which now, however, has a somewhat different interpretation.
For instance, when one has an isotropic window around the x-axis in frequency space, as in Fig.~\ref{fig:psi_3d}, top left, then this yields a disk-like window in the y-z plane  in the spatial domain, cf. Fig.~\ref{fig:psi_3d}, top right.

With Eq.~\ref{eq:psi:3d:fourier_inverse:3} we have to evaluate the radial integral to obtain a closed form expression for $\psi(x)$.
Analogous to the situation in two dimensions, for the Portilla-Simoncelli window in Eq.~\ref{eq:def:h} the integral can be computed in closed form.
The expression is similar to those in Eq.~\ref{eq:h:explicit}, see the supplementary material.

As in Sec.~\ref{sec:construction:2d}, also in three dimension the filter taps for the fast transform and other expressions like the Galerkin projection of the Laplace equation can be computed in closed form using analogous derivations.

%
%
%
%
%
%
%
%
%
%
%

\section{Proof-of-Concept Applications}
\label{sec:experiments}

\begin{figure}[t]
  \centering
  \includegraphics[width=0.48\columnwidth]{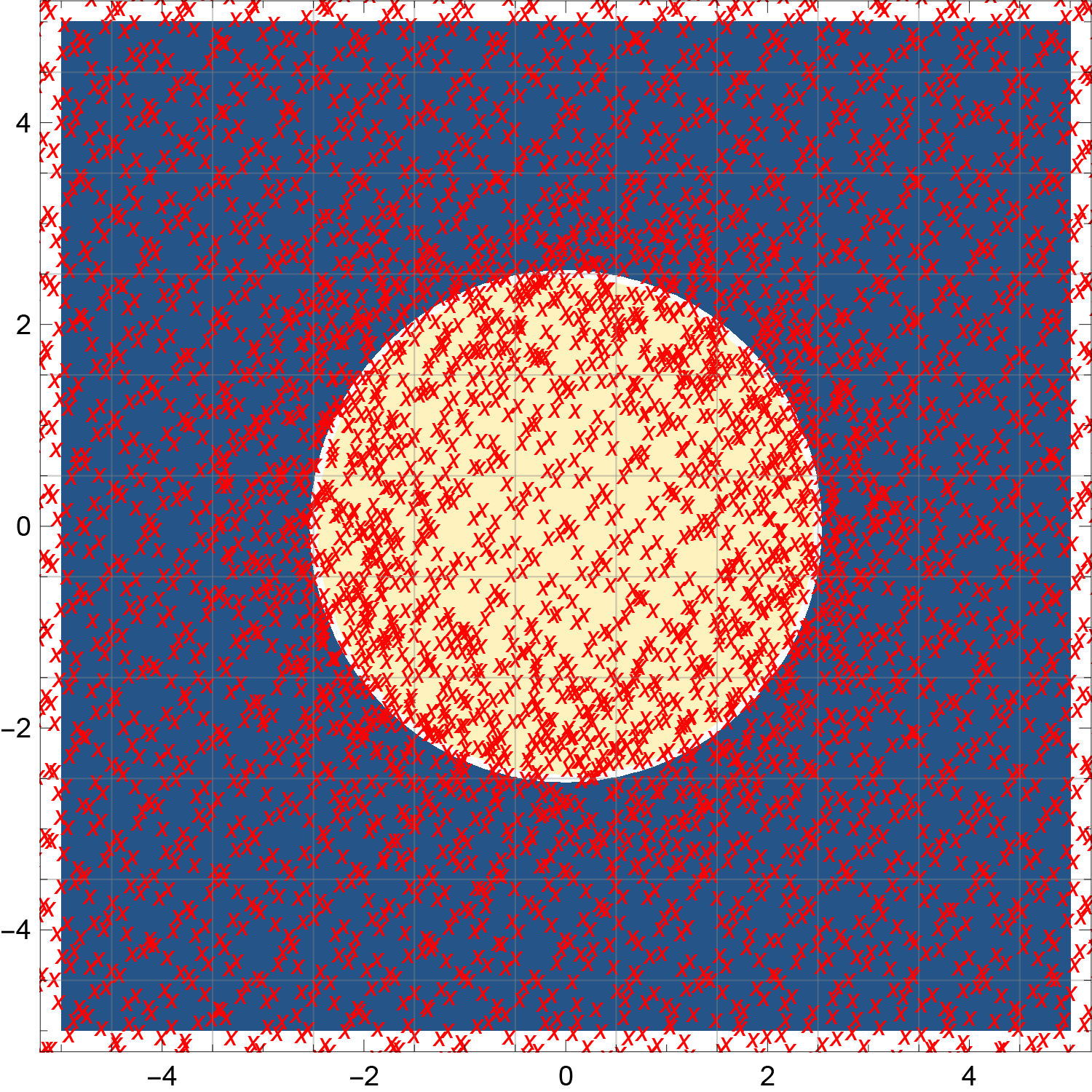}
  \includegraphics[width=0.48\columnwidth]{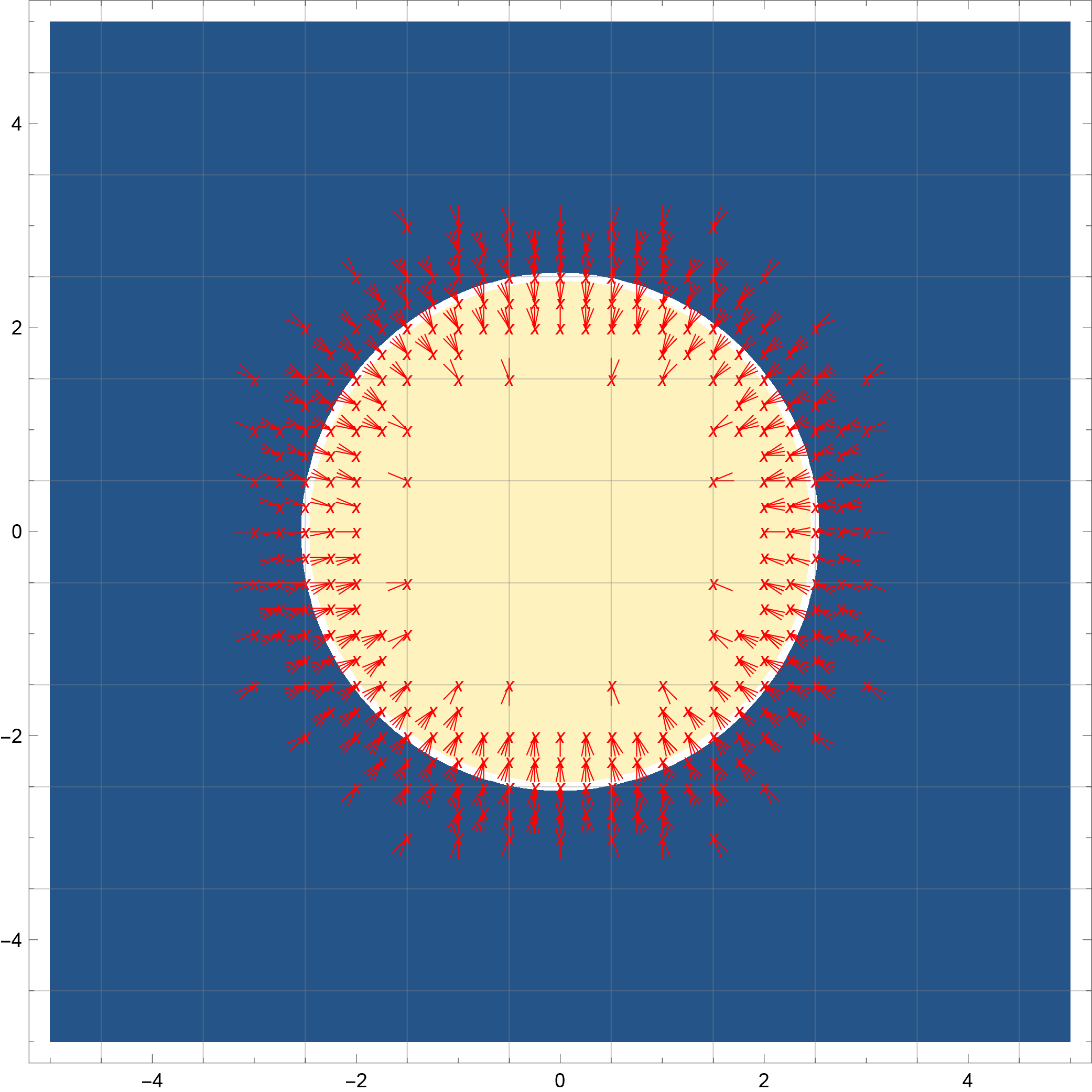}
  \caption{\emph{Left: } Adaptive set of samples. \emph{Right:} Sparse structure of basis function coefficients predicted from the vicinity to the jump discontinuity and its local orientation. As shown in the supplementary material, the sparse structure is very similar to those obtained by thresholding a full representation.}
  \label{fig:indicator:sparse_structure}
\end{figure}

\subsection{Signal estimation from point-wise samples}

Many applications require the reconstruction of signals from point-wise samples~\cite{Unser2000}.
When the data is high dimensional and/or the number of samples is very large, then the samples are often unstructured and sparse, i.e. they have a non-uniform density reflecting some of the signal's properties.
Just two examples are range scanner data, where the samples lie on a surface $\mathcal{M} \subset \mathbb{R}^3$ with a distribution induced by the geometry of $\mathcal{M}$, and samples in ray tracing-based image generation, which are often generated adaptively~\cite{Zwicker2015}, see Fig.~\ref{fig:indicator:sparse}.

In the following, we will consider the latter example and reconstruct the sparse wavelet representation of a test image from from adaptive, pointwise ray tracing samples.
We work with the wavelets from Sec.~\ref{sec:construction:2d} and use isotropic functions on coarse scales and curvelet-like ones on finer levels, as is natural for our cartoon-like test image, cf. Fig.~\ref{fig:indicator:sparse_structure}.
Following the ray tracing literature~\cite{pbrt2}, we will employ (quasi) random samples $\mathcal{X} = \{ x_i \}$ on the image plane to avoid the structured aliasing that results from uniform ones.

\begin{figure}[t]
  \centering
  \includegraphics[width=0.48\columnwidth]{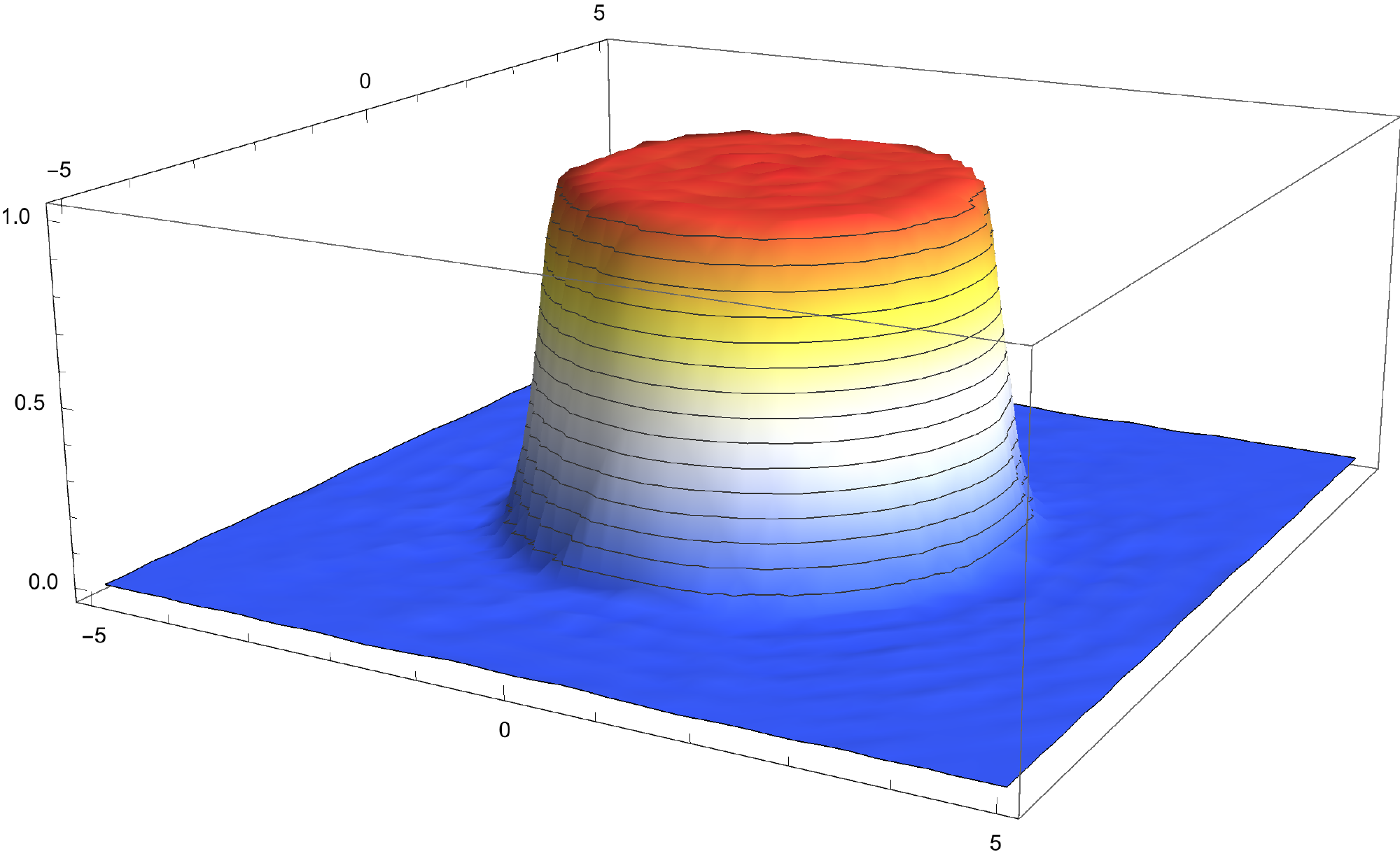}
  \includegraphics[trim={250 1100 150 0},clip,width=0.48\columnwidth]{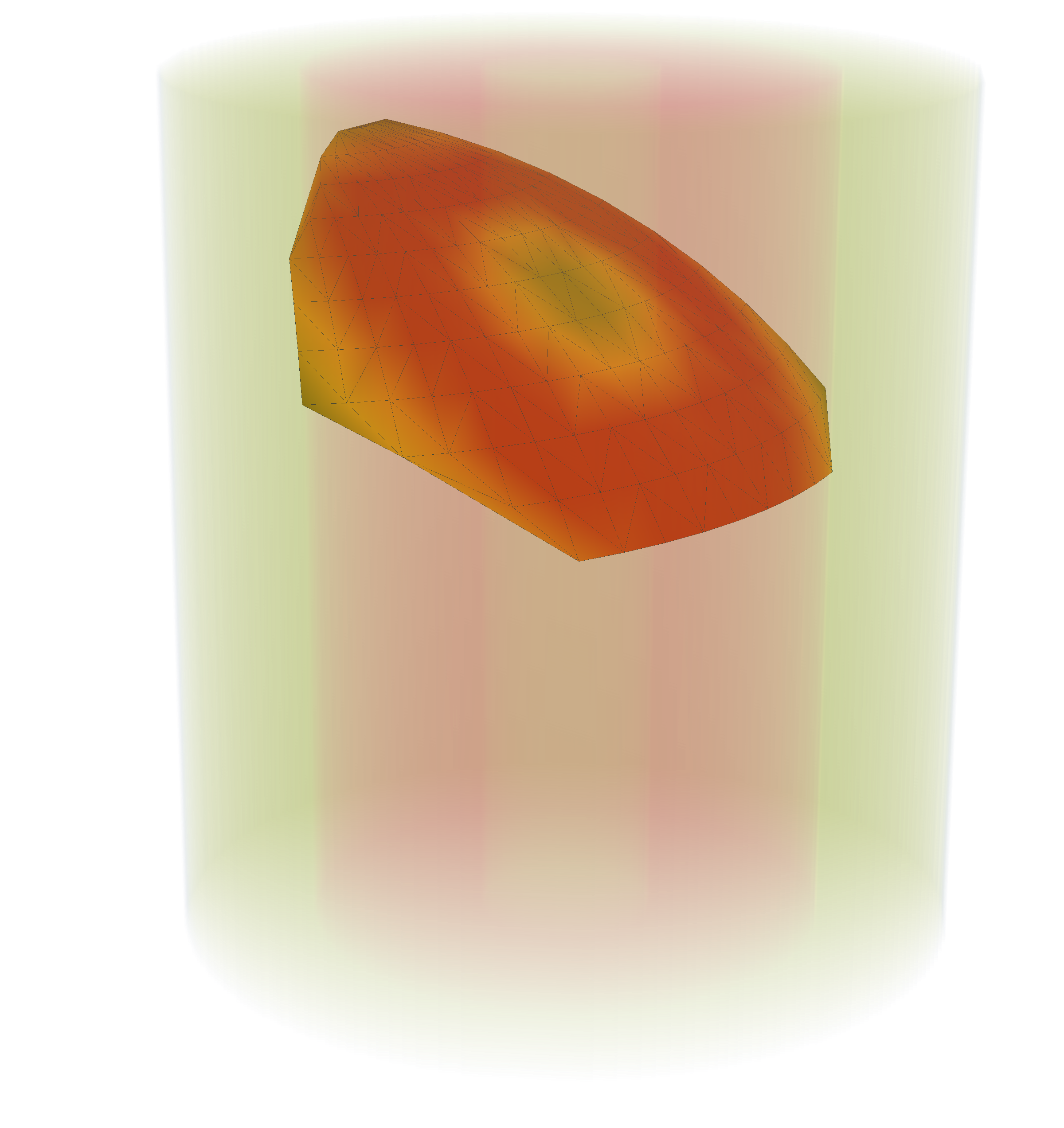}
  \caption{\emph{Left: }Reconstruction of a sparse representation of $\mathbb{B}_{r=2.5}^1$ with the wavelets of Sec.~\ref{sec:construction:2d} from 4096 sparse and adaptive ray tracing samples. Adaptive set of quasi random samples with a higher sampling density around the object boundary. \emph{Right: }Reconstruction of data on an arbitrary sub-manifold, as required for example for medical data visualization applications. The full data set is shown in the background. Such applications benefit from the closed form expressions for the frame functions in the spatial domain that enable us to evaluate them at arbitrary locations, in this case the vertices of the mesh describing the visualization sub-manifold.}
  \label{fig:indicator:sparse}
\end{figure}

The given pointwise samples $f(x_i)$ can be interpreted as expansion coefficients for a reproducing kernel representation of the image signal $f(x)$~\cite{Unser2000},
\begin{align}
  \label{eq:rk:rep}
  f(x)
  = \sum_{x_i \in \mathcal{X}} \langle f(y) , k_{x_i}(y) \rangle \, \tilde{k}_i(x)
  = \sum_{x_i \in \mathcal{X}} f(x_i) \, \tilde{k}_i(x)
\end{align}
where the $\tilde{k}_i(x)$ are dual or reconstruction kernels that are specific to the sample set $\mathcal{X}$ and the function space in which $f(x)$ lies (approximately), see~\cite{Lessig2014a} for a detailed discussion of the approach.
With Eq.~\ref{eq:rk:rep}, computing the wavelet representation of $f(x)$ amounts to solving a linear system that realizes the change of basis from the reproducing kernel frame to the wavelets,
\begin{align}
  \label{eq:rk:reconst}
  K \, f_{\psi} = f_{\mathcal{X}} ,
\end{align}
where $f_{\psi}$ is the vector of wavelet basis function coefficients and $f_{\mathcal{X}}$ those formed by the pointwise samples $f(x_i)$.
The interpolation matrix $K$ has entries $K_{si} = \psi_s(x_i)$.
Since the $x_i$ are (quasi) random, an efficient numerical implementation of Eq.~\ref{eq:rk:reconst} is facilitated by the closed form representation of the wavelets in Eq.~\ref{eq:2D:psi:hat:inverse} and Eq.~\ref{eq:h:explicit}.

After reconstrucing the wavelet coefficients using Eq.~\ref{eq:rk:reconst}, one would typically apply a form of non-linear approximation to obtain a sparse representation of $f(x)$ for storage or transmission.
However, this is inefficient because one immediately discards most of the wavelet coefficients that one just computed.
We can avoid this by using information about the location and orientation of object boundaries, which in our application is readily available through the image generation process~\cite{Zwicker2015}, and known approximation results~\cite{Candes2004,Do2005a}, which establish that only coefficients in the vicinity of high frequency features and aligned with the discontinuity will be significant.
We hence use only the corresponding basis functions and only reconstruct coefficients that are potentially non-negligible.
In other words, we work in a sparse function space that is adapted to the signal, see Fig.~\ref{fig:indicator:sparse_structure}.

Fig.~\ref{fig:indicator:sparse} shows the signal estimated from a sparse set of samples.
The test image was a (slightly regularized) indicator function of the ball $\mathbb{B}_{r=2.5}^1$.
For the experiments we used a dense set of isotropic scaling functions and isotropic wavelets on level $0$, and curvelet-like, directional ones on levels $j=1,2$, as shown in Fig.~\ref{fig:psi:2d:j1}.
For these levels there would have been $30,666$ basis functions in a dense representation while with our a priori sparse structure in Fig.~\ref{fig:indicator:sparse_structure} there are $1,332$ basis functions.
For the reconstruction we used $4,096$ adaptive samples with $3072$ being quasi uniformly distributed and $1024$ concentrated around the object boundary of the ball, see Fig.~\ref{fig:indicator:sparse_structure}.
Although some jitter can be seen in the reconstruction, the signal is overall well recovered and $L_{\infty} = 0.016$.
Since we use far fewer samples ($4,096$) than there are basis functions in the dense representation ($30,666$) the classical approach of reconstructing all coefficients and then applying nonlinear approximation would not have been possible.

\subsection{Signal reconstruction over sub-manifolds}

Applications such as curved planar reformation~\cite{Kanitsar2002} in medical imaging require the reconstruction of a signal over an arbitrary sub-manifold in $\mathbb{R}^3$.
While this can be implemented by performing a fast transform and then interpolating from the grid points closest to the sub-manifold, with the closed form spatial representations it is more efficient and straight forward to to directly reconstruct the signal on the sub-manifold.
We demonstrate this in Fig.~\ref{fig:indicator:sparse} where we visualize a bone-like dataset, consisting of nested cylinders, on a triangulated sub-manifold.
For this, we directly evaluate the basis representation of the dataset on the mesh vertices using Eq.~\ref{eq:psi:3d:fourier_inverse}.

\section{Conclusion}
\label{sec:conclusion}

In this letter we demonstrated that steerable and directional wavelets in two and three dimensions, described by Eq.~\ref{eq:2D:psi:hat} and Eq.~\ref{eq:psi:hat:3d}, respectively, have closed form expressions in the spatial domain and that also other calculations involving these, such as the filter taps for the fast transform or the Galerkin projection of the Laplace operator, can be computed analytically. 
We presented two applications that are facilitated by our closed form expressions: the estimation of a signal from pointwise samples and signal reconstruction over an arbitrary sub-manifold.
While our reconstruction procedure is very simple, e.g. compared to~\cite{Chenouard2012}, the millions of samples that one typically uses on the image plane in ray tracing also require a very efficient method.
Our application also demonstrates a principle that we believe will be of considerable importance for high dimensional data in the future: the reconstruction of a sparse signal representation from a sparse set of samples. 

Our closed form expressions for the radial window in the spatial domain are rather complicated.
In future work we would like to find windows with simpler expressions and potentially also with faster decay in space, for example following the approach proposed in~\cite{Ward2015}.

A reference implementation for the work discussed in this letter is available in the supplementary material.

\section*{Acknowledgment}
The author would like to thank Eugene Fiume and Marc Alexa for continuing support.

\ifCLASSOPTIONcaptionsoff
  \newpage
\fi



\bibliographystyle{IEEEtran}
\bibliography{polarlets}

\begin{thebibliography}{10}
\providecommand{\url}[1]{#1}
\csname url@samestyle\endcsname
\providecommand{\newblock}{\relax}
\providecommand{\bibinfo}[2]{#2}
\providecommand{\BIBentrySTDinterwordspacing}{\spaceskip=0pt\relax}
\providecommand{\BIBentryALTinterwordstretchfactor}{4}
\providecommand{\BIBentryALTinterwordspacing}{\spaceskip=\fontdimen2\font plus
\BIBentryALTinterwordstretchfactor\fontdimen3\font minus
  \fontdimen4\font\relax}
\providecommand{\BIBforeignlanguage}[2]{{%
\expandafter\ifx\csname l@#1\endcsname\relax
\typeout{** WARNING: IEEEtran.bst: No hyphenation pattern has been}%
\typeout{** loaded for the language `#1'. Using the pattern for}%
\typeout{** the default language instead.}%
\else
\language=\csname l@#1\endcsname
\fi
#2}}
\providecommand{\BIBdecl}{\relax}
\BIBdecl

\bibitem{Freeman1991}
\BIBentryALTinterwordspacing
W.~T. Freeman and E.~H. Adelson, ``{The design and use of steerable filters},''
  \emph{IEEE Transactions on Pattern Analysis and Machine Intelligence},
  vol.~13, no.~9, pp. 891--906, 1991. [Online]. Available:
  \url{http://ieeexplore.ieee.org/document/93808/}
\BIBentrySTDinterwordspacing

\bibitem{Perona1991}
\BIBentryALTinterwordspacing
P.~Perona, ``{Deformable kernels for early vision},'' in \emph{Proceedings.
  1991 IEEE Computer Society Conference on Computer Vision and Pattern
  Recognition}.\hskip 1em plus 0.5em minus 0.4em\relax IEEE Comput. Sco. Press,
  1991, pp. 222--227. [Online]. Available:
  \url{http://ieeexplore.ieee.org/document/139691/}
\BIBentrySTDinterwordspacing

\bibitem{Unser2011}
\BIBentryALTinterwordspacing
M.~Unser, N.~Chenouard, and D.~{Van De Ville}, ``{Steerable Pyramids and Tight
  Wavelet Frames in {\$}L{\_}{\{}2{\}}({\{}R{\}}{\^{}}{\{}d{\}}){\$}},''
  \emph{IEEE Transactions on Image Processing}, vol.~20, no.~10, pp.
  2705--2721, oct 2011. [Online]. Available:
  \url{http://ieeexplore.ieee.org/document/5746534/}
\BIBentrySTDinterwordspacing

\bibitem{Chenouard2012}
\BIBentryALTinterwordspacing
N.~Chenouard and M.~Unser, ``{3D Steerable Wavelets in Practice},'' \emph{IEEE
  Transactions on Image Processing}, vol.~21, no.~11, pp. 4522--4533, nov 2012.
  [Online]. Available: \url{http://ieeexplore.ieee.org/document/6226458/}
\BIBentrySTDinterwordspacing

\bibitem{Unser2013}
\BIBentryALTinterwordspacing
M.~Unser and N.~Chenouard, ``{A Unifying Parametric Framework for 2D Steerable
  Wavelet Transforms},'' \emph{SIAM Journal on Imaging Sciences}, vol.~6,
  no.~1, pp. 102--135, jan 2013. [Online]. Available:
  \url{http://epubs.siam.org/doi/10.1137/120866014}
\BIBentrySTDinterwordspacing

\bibitem{Ward2014}
\BIBentryALTinterwordspacing
J.~P. Ward and M.~Unser, ``{Harmonic singular integrals and steerable wavelets
  in L2(Rd)},'' \emph{Applied and Computational Harmonic Analysis}, vol.~36,
  no.~2, pp. 183--197, mar 2014. [Online]. Available:
  \url{https://www.sciencedirect.com/science/article/pii/S1063520313000274}
\BIBentrySTDinterwordspacing

\bibitem{Candes2005a}
\BIBentryALTinterwordspacing
E.~J. Cand{\`{e}}s and D.~L. Donoho, ``{Continuous curvelet transform: I.
  Resolution of the Wavefront Set},'' \emph{Applied and Computational Harmonic
  Analysis}, vol.~19, no.~2, pp. 162--197, sep 2005. [Online]. Available:
  \url{http://www.sciencedirect.com/science/article/pii/S1063520305000199}
\BIBentrySTDinterwordspacing

\bibitem{Candes2005b}
\BIBentryALTinterwordspacing
------, ``{Continuous curvelet transform: II. Discretization and Frames},''
  \emph{Applied and Computational Harmonic Analysis}, vol.~19, no.~2, pp.
  198--222, sep 2005. [Online]. Available:
  \url{http://www.sciencedirect.com/science/article/pii/S1063520305000205}
\BIBentrySTDinterwordspacing

\bibitem{Do2005a}
\BIBentryALTinterwordspacing
M.~Do and M.~Vetterli, ``\BIBforeignlanguage{English}{{The contourlet
  transform: an efficient directional multiresolution image representation}},''
  \emph{\BIBforeignlanguage{English}{IEEE Transactions on Image Processing}},
  vol.~14, no.~12, pp. 2091--2106, dec 2005. [Online]. Available:
  \url{http://ieeexplore.ieee.org/articleDetails.jsp?arnumber=1532309}
\BIBentrySTDinterwordspacing

\bibitem{Labate2005}
\BIBentryALTinterwordspacing
D.~Labate, W.-Q. Lim, G.~Kutyniok, and G.~Weiss, ``{Sparse Multidimensional
  Representation using Shearlets},'' in \emph{Wavelets XI}, M.~Papadakis, A.~F.
  Laine, and M.~A. Unser, Eds.\hskip 1em plus 0.5em minus 0.4em\relax
  International Society for Optics and Photonics, aug 2005, pp. 254--262.
  [Online]. Available:
  \url{http://proceedings.spiedigitallibrary.org/proceeding.aspx?articleid=870627}
\BIBentrySTDinterwordspacing

\bibitem{Portilla2000}
\BIBentryALTinterwordspacing
J.~Portilla and E.~P. Simoncelli, ``{A Parametric Texture Model Based on Joint
  Statistics of Complex Wavelet Coefficients},'' \emph{International Journal of
  Computer Vision}, vol.~40, no.~1, pp. 49--70, 2000. [Online]. Available:
  \url{http://link.springer.com/10.1023/A:1026553619983}
\BIBentrySTDinterwordspacing

\bibitem{Papadakis2003}
\BIBentryALTinterwordspacing
M.~Papadakis, G.~Gogoshin, I.~A. Kakadiaris, D.~J. Kouri, and D.~K. Hoffman,
  ``{Nonseparable Radial Frame Multiresolution Analysis in Multidimensions},''
  \emph{Numerical Functional Analysis and Optimization}, vol.~24, no. 7-8, pp.
  907--928, dec 2003. [Online]. Available:
  \url{http://www.tandfonline.com/doi/abs/10.1081/NFA-120026385}
\BIBentrySTDinterwordspacing

\bibitem{Azencott2009}
\BIBentryALTinterwordspacing
R.~Azencott, B.~G. Bodmann, and M.~Papadakis, ``{Steerlets: a novel approach to
  rigid-motion covariant multiscale transforms},'' V.~K. Goyal, M.~Papadakis,
  and D.~{Van De Ville}, Eds., aug 2009, p. 74460A. [Online]. Available:
  \url{http://proceedings.spiedigitallibrary.org/proceeding.aspx?doi=10.1117/12.825704}
\BIBentrySTDinterwordspacing

\bibitem{Romero2009a}
\BIBentryALTinterwordspacing
J.~R. Romero, S.~K. Alexander, S.~Baid, S.~Jain, and M.~Papadakis, ``{The
  geometry and the analytic properties of isotropic multiresolution
  analysis},'' \emph{Advances in Computational Mathematics}, vol.~31, no. 1-3,
  pp. 283--328, oct 2009. [Online]. Available:
  \url{http://link.springer.com/10.1007/s10444-008-9111-6}
\BIBentrySTDinterwordspacing

\bibitem{Walter1999}
\BIBentryALTinterwordspacing
G.~G. Walter and L.~Cai, ``{Periodic Wavelets from Scratch},'' \emph{Journal of
  Computational Analysis and Applications}, vol.~1, no.~1, pp. 25--41, 1999.
  [Online]. Available: \url{http://link.springer.com/10.1023/A:1022614519335}
\BIBentrySTDinterwordspacing

\bibitem{Narcowich2006a}
\BIBentryALTinterwordspacing
F.~J. Narcowich, P.~Petrushev, and J.~D. Ward, ``{Localized Tight Frames on
  Spheres},'' \emph{SIAM Journal on Mathematical Analysis}, vol.~38, no.~2, pp.
  574--594, jan 2006. [Online]. Available:
  \url{http://epubs.siam.org/doi/10.1137/040614359}
\BIBentrySTDinterwordspacing

\bibitem{Freeden1998}
W.~Freeden, T.~Gervens, and M.~Schreiner, \emph{{Constructive Approximation on
  the Sphere (With Applications to Geomathematics)}}.\hskip 1em plus 0.5em
  minus 0.4em\relax Oxford Sciences Publication. Clarendon Press, Oxford
  University, 1998.

\bibitem{McEwen2016}
\BIBentryALTinterwordspacing
J.~D. McEwen, C.~Durastanti, and Y.~Wiaux, ``{Localisation of directional
  scale-discretised wavelets on the sphere},'' \emph{Applied and Computational
  Harmonic Analysis}, 2016. [Online]. Available:
  \url{http://www.sciencedirect.com/science/article/pii/S1063520316000324}
\BIBentrySTDinterwordspacing

\bibitem{Kutyniok2016}
G.~Kutyniok and P.~Petersen, ``{Classification of Edges Using Compactly
  Supported Shearlets},'' \emph{Applied and Computational Harmonic Analysis},
  nov 2016.

\bibitem{Unser2000}
M.~Unser, ``{Sampling---50 Years After Shannon},'' \emph{Proceedings of the
  IEEE}, vol.~88, pp. 569--587, apr 2000.

\bibitem{Zwicker2015}
\BIBentryALTinterwordspacing
M.~Zwicker, W.~Jarosz, J.~Lehtinen, B.~Moon, R.~Ramamoorthi, F.~Rousselle,
  P.~Sen, C.~Soler, and S.-E. Yoon, ``{Recent Advances in Adaptive Sampling and
  Reconstruction for Monte Carlo Rendering},'' \emph{Computer Graphics Forum},
  vol.~34, no.~2, pp. 667--681, may 2015. [Online]. Available:
  \url{http://doi.wiley.com/10.1111/cgf.12592}
\BIBentrySTDinterwordspacing

\bibitem{pbrt2}
M.~Pharr and G.~Humphreys, \emph{{Physically Based Rendering: From Theory to
  Implementation}}, 2nd~ed.\hskip 1em plus 0.5em minus 0.4em\relax San
  Francisco, CA, USA: Morgan Kaufmann Publishers Inc., 2010.

\bibitem{Lessig2014a}
\BIBentryALTinterwordspacing
C.~Lessig, M.~Desbrun, and E.~Fiume, ``{A Constructive Theory of Sampling for
  Image Synthesis Using Reproducing Kernel Bases},'' \emph{ACM Transactions on
  Graphics (Proceedings of SIGGRAPH 2014)}, vol.~33, no.~4, pp. 1--14, jul
  2014. [Online]. Available:
  \url{http://dl.acm.org/citation.cfm?id=2601097.2601149}
\BIBentrySTDinterwordspacing

\bibitem{Candes2004}
\BIBentryALTinterwordspacing
E.~Cand{\`{e}}s and D.~L. Donoho, ``{New Tight Frames of Curvelets and Optimal
  Representations of Objects with Piecewise {\$}C{\^{}}2{\$} Singularities},''
  \emph{Communications on Pure and Applied Mathematics}, vol.~57, no.~2, pp.
  219--266, feb 2004. [Online]. Available:
  \url{http://doi.wiley.com/10.1002/cpa.10116}
\BIBentrySTDinterwordspacing

\bibitem{Kanitsar2002}
\BIBentryALTinterwordspacing
A.~Kanitsar, D.~Fleischmann, R.~Wegenkittl, P.~Felkel, and E.~Groller, ``{CPR -
  curved planar reformation},'' in \emph{IEEE Visualization, 2002. VIS
  2002.}\hskip 1em plus 0.5em minus 0.4em\relax IEEE, 2002, pp. 37--44.
  [Online]. Available: \url{http://ieeexplore.ieee.org/document/1183754/}
\BIBentrySTDinterwordspacing

\bibitem{Ward2015}
\BIBentryALTinterwordspacing
J.~P. Ward, P.~Pad, and M.~Unser, ``{Optimal Isotropic Wavelets for Localized
  Tight Frame Representations},'' \emph{IEEE Signal Processing Letters},
  vol.~22, no.~11, pp. 1918--1921, nov 2015. [Online]. Available:
  \url{http://ieeexplore.ieee.org/document/7130607/}
\BIBentrySTDinterwordspacing

\end{thebibliography}
\end{document}